\documentclass[12pt]{article}

\usepackage{amssymb}
\usepackage{amsmath}
\usepackage{amsbsy}
\usepackage{amscd}
\usepackage{amsfonts}
\usepackage{amsthm}
\usepackage{mathrsfs}
\usepackage{verbatim}
\usepackage[colorlinks]{hyperref}
\usepackage{fullpage}
\usepackage{mathdots}
\usepackage{graphicx,subfigure}
\usepackage[english]{babel}
\usepackage[utf8]{inputenc}
\usepackage{tikz}
\usepackage{adjustbox}
\usepackage{authblk}
\usepackage{caption}

\usepackage{mathtext}

\usepackage{tikz-cd}
\usetikzlibrary{backgrounds,fit, matrix}
\usetikzlibrary{positioning}
\usetikzlibrary{calc,through,chains}
\usetikzlibrary{arrows,shapes,snakes,automata, petri}
\usepackage[boxsize=1.25em, centerboxes]{ytableau}

\newtheorem{Theorem}[subsection]{Theorem}
\newtheorem{Corollary}[subsection]{Corollary}
\newtheorem{Lemma}[subsection]{Lemma}

\theoremstyle{definition}
\newtheorem{Definition}[subsection]{Definition}

\newtheorem{Question}[subsection]{Question}
\newtheorem{Remark}[subsection]{Remark}

\numberwithin{equation}{section}

\newcommand{\PP}{{\mathbb P}}
\newcommand{\C}{{\mathbb C}}

\newcommand{\N}{{\mathbb N}}

\DeclareMathOperator{\Pic}{Pic}

\DeclareMathOperator{\GL}{\textbf{GL}}
\DeclareMathOperator{\SL}{\textbf{SL}}
\DeclareMathOperator{\Sp}{\textbf{Sp}}
\DeclareMathOperator{\SO}{\textbf{SO}}
\DeclareMathOperator{\Spin}{\textbf{Spin}}

\newcommand{\X}{\mathcal{X}}
\newcommand{\U}{\mathscr{U}}

\begin{document}

\title{On the Fixed Points of a Regular Unipotent Element}

\author[1]{Mahir Bilen Can}

\affil[1]{\normalsize{
Tulane University, New Orleans, Louisiana\\
mahirbilencan@gmail.com}} 

\maketitle

\begin{abstract} 
The fixed point variety of a regular unipotent element on a wonderful completion is investigated. 
For the wonderful completion of the quotient by a symmetric Levi subgroup, it is shown that the fixed point variety is $SL_2$-regular. 
\\



\end{abstract}

\section{Introduction}

Let $G$ be a semi-simple simply connected complex algebraic group. 
It is a remarkable theorem of Springer~\cite{Springer1969} that the following natural map is a resolution of singularities:
\begin{align}
pr_1 : \{ (u,gB) \in \U(G) \times G/B:\ ugB = gB \} &\longrightarrow \U(G)\\
(u,gB)& \longmapsto u \notag,
\end{align}
where $\U(G)$ is the variety of unipotent elements in $G$, and $B$ is a Borel subgroup of $G$. 
Then the singularity of a point $u$ in $\U(G)$ is measured by the fiber $pr_1^{-1}(u)$, which is called the {\em Springer fiber at $u$}. 
The top dimensional cohomology space of the Springer fiber affords a representation of the Weyl group. 
In this note, we investigate a special case of a family of fibers that are defined similarly to the Springer fibers in a closely related setup. 
More precisely, we consider the fixed point variety of the regular unipotent element of $G$, which acts on a complex Hermitian  symmetric space of noncompact type. Our main goal is to prove the following statement. 
\begin{Theorem}\label{T:firstthm}
Let $X$ denote the wonderful completion a complex Hermitian symmetric space of the form $G/N_G(K)$, where $N_G(K)$ is the normalizer of the fixed subgroup of an automorphism of $G$ of period 2. 
Then there exists a regular $\SL_2$ action on $X$.
\end{Theorem}
In this theorem, the term {\em regular} means that a maximal unipotent subgroup of $\SL_2$ has an isolated fixed point in $X$. 
As shown in~\cite{AC1987}, for such an action of $\SL_2$ on a smooth projective variety $X$, if $Z$ denotes the zero scheme of the vector field induced from the action of the maximal unipotent subgroup of $\SL_2$, then the coordinate ring $A(Z)$ has a natural grading, making it isomorphic as a graded ring to $H^*(X,\C)$. Furthermore, there is a remarkable equivariant extension of this beautiful fact.
In~\cite{BrionCarrell}, Brion and Carrell show that the equivariant cohomology ring of an $\SL_2$-regular variety $X$,
which is projective and smooth, is isomorphic to the coordinate ring of an affine curve in $X\times \PP^1$.  
We should mention that, in general, wonderful completions are not $\SL_2$-regular varieties~\cite{CanJoyce,CarrellKiumars}
since the fixed point set of a regular unipotent element can be positive dimensional. 
In fact, the Poincar\'e polynomial of the $u$-fixed subvariety of $\overline{\SL_n/N_{\SL_n}(\SO_n)}$ 
(resp. of $\overline{\SL_n/N_{\SL_n}(\Sp_n)}$) is given by $\sum_{i=0}^n {n-i-1 \choose i} t^{2i}$ (resp. by $(1+t^2)^{n/2-1}$), see~\cite{CHJ2}. 
However, in none of these examples the semi-simple rank of the subgroup $N_G(K)$ is equal to the semi-simple rank of $G$; this plays a role in the regularity of the action but it does not completely characterize it.

To prove our main result (Theorem~\ref{T:firstthm}) we prove an auxiliary theorem.   
Let $u$ be a unipotent element of $G$, and let $Y$ be a $G$-variety. 
We denote by $Y_u$ the possibly reduced subvariety of $u$-fixed points in $Y$. 
It is well-known that~\cite{Steinberg1974} if $Y$ is a homogeneous space of the form $G/P$, where $P$ is a parabolic subgroup of $G$, then $Y_u$ is a singleton whenever $u$ is regular. Let $L$ be a Levi subgroup of $P$. 
Our second main result is the following. 
\begin{Theorem}\label{T:secondthm}
If $u$ is a regular unipotent element of $G$, then $u$ acts freely on $G/J$, where $J$ is any closed subgroup of the normalizer of $L$ in $G$. 
\end{Theorem}
This result should be viewed as a first step towards studying $\SL_2$-regular actions on the equivariant embeddings of a large class of  homogeneous spaces. In the sequel, we will explain this comment in more detail.

One of our main motivations for proving Theorem~\ref{T:firstthm} is that the cohomology rings of wonderful completions of symmetric spaces have rather intriguing but complicated combinatorial structures. 
Their Poincar\'e polynomials are computed by De Concini and Springer in~\cite{DS} as (highly nontrivial) sums of polynomials
with contributions from the orbits that contain torus fixed points. 
We want to simplify their summation formulas; if it is possible, we want to express them as products. 
Of course, it is unreasonable to expect a simple factorization formula for every integral polynomial.
Nevertheless, as shown by Aky\i ld\i z and Carrell in~\cite{AC1989}, the Poincar\'e polynomials of $\SL_2$-regular varieties have nice factorization properties. 
Indeed, if $X$ is a smooth projective $\SL_2$-regular variety, then its Poincar\'e polynomial is given by 
\begin{align*}
P(X,t) = \prod_{1\leq i \leq k } \left( \frac{1-t^{2a_i+2}}{1-t^{2a_i}} \right)^{\mu_i}, 
\end{align*}
where $\mu_i$'s ($i\in \{1,\dots, k\}$) are the dimensions of the $\C^*$-weight subspaces of the tangent space at the unique fixed point of 
the maximal unipotent subgroup of $\SL_2$, and the $a_i$'s ($i\in \{1,\dots, k\}$) are the negatives of the corresponding $\C^*$-weights. 
In other words, roughly speaking, the information of the cohomology ring of a smooth projective $\SL_2$-variety is encoded at the 
zero scheme of the unique fixed point of the maximal unipotent subgroup of $\SL_2$. 
Now the following result should not come as a big surprise.
\begin{Corollary}\label{C:firstcorollary}
Let $X$ be the wonderful completion of the symmetric space $G/L$, where $L$ is a Levi subgroup of a parabolic subgroup of $G$. 
Let $Y$ denote the unique closed $G$-orbit in $X$.  Then the Poincar\'e polynomial of $X$ is divisible by the Poincar\'e polynomial of $Y$. 
In fact, we have $P(X,t) = P(Y,t)\prod_{ i=1}^r \frac{1-t^{2(i+1)}}{1-t^{2i}}$, where $r$ is the number of prime $G$-stable divisors in $X$. 
\end{Corollary}
This corollary is closely related to the fine structure of the torus weights on the normal space at the unique fixed point of the maximal unipotent subgroup. 
We will present the full scope of this calculation somewhere else.

Before we describe the organization of our manuscript, we want to make a final remark about the extendability of our first main result. 
A closed subgroup $H\subseteq G$ is called a {\em spherical subgroup} if under the left multiplication action $G\times G/H \to G/H$, $(g,aH)\mapsto gaH$,  
a Borel subgroup of $G$ has an open orbit in $G/H$. 
More generally, a normal irreducible $G$-variety $X$ is said to be a {\em spherical $G$-variety} if a Borel subgroup of $G$ has an open orbit in $X$. 
It is well-known that the symmetric subgroups are spherical~\cite{Matsuki1979} (over $\C$)~\cite{Springer1985} (in general).
There is a complete classification of the symmetric spaces~\cite{Springer1987} with origins in \'E. Cartan's work.
In fact, for semi-simple $G$, the complete classification of affine spherical homogeneous spaces is well-known.
This classification has been achieved by Kr\"amer~\cite{Kramer1979} for simple $G$, and by Brion~\cite{Brion1987} and Mikityuk~\cite{Mikityuk1986} for non-simple semi-simple $G$. 
It is shown in~\cite{Brion1987} that if $G/H$ is an affine spherical variety, then $H$ is contained in a Levi subgroup of $G$. 
Therefore, for every reductive spherical subgroup $H$ of $G$, our Theorem~\ref{T:secondthm} is applicable. 
Notice that, since every closed orbit in an equivariant embedding contain a $u$-fixed point, our theorem is a necessary but not sufficient condition for the $\SL_2$-regularity of a complete equivariant embedding of $G/H$. 
It would be interesting to classify all (reductive) spherical subgroups $H$ such that $G/H$ has an $\SL_2$-regular wonderful completion.

The organization of this note is as follows. 
In the next section, we setup our notation.
Also, we review some standard facts about wonderful completions, spherical affine homogeneous spaces,
and $\SL_2$-regular varieties. 
We prove our second theorem in Section~\ref{S:second}, and we prove our first main result as well as its Corollary~\ref{C:firstcorollary} in Section~\ref{S:first}.

\section{Preliminaries}

Although most of our statements hold true over algebraically closed fields of arbitrary characteristics, for simplicity, 
we assume that all of our algebraic varieties and schemes are defined over $\C$. 
The notation $\mathbb{G}_a$ (resp. $\mathbb{G}_m$) will be used for the abelian group $(\C,+)$ (resp. $(\C^*,\cdot)$). 
The unipotent radical of an algebraic group will be denoted by $\mathscr{R}_uG$. 
Let $T$ be a maximal torus in a reductive group $G$. 
The Weyl group of $(G,T)$, that is, $N_G(T)/T$, will be denoted by $W$. 
If $L$ is a Levi subgroup of $G$, then its Weyl group is a subgroup of $W$; it will be denoted by $W_L$ (or by $W_I$, where $I$
is a set of simple roots that determine $L$).
The Weyl groups are equipped with a length function, $\ell : W \to \N$, defined as follows. 
Let $B$ be a Borel subgroup of $G$ such that $T\subseteq B$. 
If $w$ is an element of $W$, then $\ell(w)$ is given by the dimension of the $B$-orbit $B\dot{w}B$ in $G/B$, where $\dot{w}$ is any representative for $w$ 
in $N_G(T)$.

\subsection{Spherical Levi subgroups.}\label{SS:reductivespherical}

We begin with some general remarks. 
Let $G$ be a connected reductive group. 
The classification of affine spherical homogeneous spaces $G/H$ easily reduces to the classification of 
the affine spherical homogeneous spaces of the form $G'/G'\cap H$, where $G'$ is the semi-simple commutator subgroup of $G$. 
Let us proceed with the assumption that $G$ is semi-simple, and $K$ is a spherical reductive subgroup of $G$. 
Then, as we mentioned earlier, the homogeneous spaces $G/K$ are classified by Brion and Mikityuk. 
Let $\text{rank}(K)$ denote the dimension of a maximal torus in $K$. 
In~\cite{Brion1987}, Brion shows that if $L$ is a spherical Levi subgroup of $G$, where $G$ is simply connected, then the couple $(G,L)$ is a product of the couples $(G',L')$ that appear in the following table:

\begin{table}[ht]
\begin{center}
\begin{tabular}{|l |c|c| }
  \hline
&  $G'$ & $L'$ \\ 
  \hline
1.&  $\SL_{p+q}$ & $\SL_{p+q} \cap (\GL_p \times \GL_q)$ \\
  \hline
2.&  $\Sp_{2n}$ & $\mathbb{G}_m \times \Sp_{2n-2}$ \\
  \hline 
3.&  $\Sp_{2n}$ & $\GL_n$ \\
  \hline
4.&  $\Spin_{n},\ n\geq 5$ & $\mathbb{G}_m \times \SO_{n-2}$ \\
  \hline
5.&  $\Spin_{2n}$ & $\GL_n$ \\
  \hline
6.&  $\Spin_{2n+1}$ & $\GL_n$ \\
  \hline
7.&  $E_6$ & $\mathbb{G}_m \times D_5$ \\
  \hline
8.&  $E_7$ & $\mathbb{G}_m \times E_6$ \\
  \hline
\end{tabular}
\end{center}
\caption{Irreducible spherical Levi subgroups.}
\label{T:Classical}
\end{table}

In~\cite{Brundan1998}, Brundan showed that the above list of spherical Levi subgroups remains unchanged in positive characteristic. 
We note also that all of the Levi subgroups $L'$ that appear in Table~\ref{T:Classical} are reductive with a one dimensional connected center. 
This information is useful for the computation of the Picard group of the homogeneous space $G/L$ as well as its wonderful completion.

\subsection{Symmetric spaces of Hermitian type.}\label{SS:Hermitian}

A homogeneous space of the form $G/K$, where $K$ is the fixed subgroup of an algebraic group automorphism $\sigma : G\to G$ such that $\sigma^2 = id$ is called a {\em symmetric space}. 
In this case, $K$ will be called as a {\em symmetric subgroup of $G$}. 
A {\em complex Hermitian symmetric space of non-compact type} is a symmetric space $G/K$, 
where $G$ is a semi-simple group, and $K$ is a Levi subgroup of a parabolic subgroup in $G$. 
Such a symmetric space has a decomposition of the form $G_1 / K_1\times \cdots \times G_r/K_r$, 
where $K_i := K\cap G_i$ ($i\in \{1,\dots, r\}$), and $G_i$ ($i\in \{1,\dots, r\}$) is a (normal) simple factor of $G$.

\begin{Lemma}
Let $G$ be a simple algebraic group, and let $L$ be a symmetric subgroup of $G$. 
If $L$ is a Levi subgroup of a parabolic subgroup in $G$, then $G/L$ is one of the symmetric spaces in Table~\ref{T:Classical}.
\end{Lemma}
\begin{proof}
For simple algebraic groups over algebraically closed fields of arbitrary characteristic ($\neq 2$), the classification of involutions,
hence of the symmetric spaces, is well-known~\cite{Araki1962,Satake1971,Springer1987, Helminck1988}. 
A detailed table of the possibilities is presented in Timashev's book~\cite[Table 26.3]{Timashev}.
In this table, we search for $G/L$, where $L$ is a Levi subgroup of a parabolic subgroup of $G$. 
It is easily verified that Brion's list of spherical Levi subgroups from Table~\ref{T:Classical} is precisely the list of symmetric spaces $G/L$,
where $L$ is a Levi subgroup. 
\end{proof}


\subsection{Wonderful Completions.}\label{SS:wonderful embeddings}
		
In this subsection, we will briefly review the theory of wonderful embeddings.
We recommend the article of Pezzini~\cite{Pezzini2018} for further details.
\begin{Definition}\label{D:wonderful}
Let $X$ be an irreducible, smooth, complete $G$-variety.
If the following properties hold, then $X$ is called a {\em wonderful variety of rank $r$}: 
\begin{enumerate}
\item There is an open $G$-orbit $X_0$ whose complement is the union of smooth prime $G$-divisors $D_1,\dots, D_r$ with normal crossings.
\item The intersection $D_1\cap \cdots \cap D_r$ is nonempty.
\item If $x$ and $x'$ are such that $\{i:\ x\in D_i \} = \{i :\ x'\in D_i\}$, then $G\cdot x = G\cdot x'$. 
\end{enumerate}
\end{Definition}
If $H$ is a closed subgroup such that $G/H \cong X_0$, then $H$ is called a {\em wonderful (spherical) subgroup}.
In this case, $X$ is called the {\em wonderful completion of $G/H$}; it is unique up to isomorphism. 
If $H$ is an arbitrary spherical subgroup of $G$, then it may not be a wonderful subgroup.
Nevertheless, the normalizer $N_G(H)$ of a spherical subgroup is always wonderful by a result of Knop~\cite{Knop1996}.

Let $X$ be a wonderful variety of rank $r$. Let $Y$ denote its closed orbit.
Then $Y$ is a partial flag variety of the form $G/Q$, where $Q$ is a parabolic subgroup.
The stabilizer of any point in $Y$ is conjugate-isomorphic to $Q$. Let $P$ be a parabolic subgroup that is opposite to $Q$. 
Let $B$ be a Borel subgroup contained in $P$, and let $B^-$ denote the opposite Borel subgroup. 
Then there is a unique $B^-$-fixed point $y$ in $Y$.
Let $T$ be the maximal torus defined by $T := B\cap B^-$. 
The $T$-weights at the normal space $T_y X / T_y Y$ are called the {\em spherical roots of $X$}. 
We denote the set of spherical roots (relative to the choices that we have made) by $\varSigma (X)$. 
The set of all subsets of $\varSigma(X)$ is in order reversing bijection with the set of $G$-orbit closures in $X$. 
The Poincar\'e polynomial $P(X,t)$ of $X$ is calculated in~\cite[Theorem 2.4]{DS}. 
It is given by a formula which is rather heavy in notation.
We copy their formula here, without explaining much of its details, to make the point that $P(X,t)$ cannot be readily factorized in a geometrically meaningful way: 
$P(X,t) = \sum_{\Gamma \subseteq \varSigma(X)} t^{ r - |\Gamma| + |R_{\Omega(\Gamma)}^+ \cap \Phi_1| 
\frac{|W_{\Omega(\Gamma)}|}{m_\Gamma}} \sum_{w\in W/W_{\Omega(\Gamma)}} t^{2(\ell(w) + s_\Gamma(w))}$.
In this formula, the first summation is taken over the subsets for which the corresponding $G$-orbit contains a torus fixed point. 
Of course, if all torus fixed points are contained in the unique closed orbit, then the formula is simplified. 
This is a rare exception but it happens if $X$ is the wonderful completion of the adjoint form of $G$.
Then we have $P(X,t) = (\sum_{w\in W} t^{2\ell(w)}) (\sum_{w\in W} t^{2(\ell(w) + |Des_L(w)|)})$, where $Des_L(w)$ is the set of simple reflections such that 
$\ell(s_\alpha w)  < \ell(w)$. 
This formula for the Poincar\'e polynomial of $X = \overline{G_{ad}}$ is attributed to Lusztig in~\cite{DP}.
By using different methods, in~\cite[Theorem 4.1]{Renner2003}, Renner extended Lusztig's formula to all $G$-orbit closures. 
Let us also mention that, in~\cite{CarrellKiumars}, Carrell and Kaveh show that $X = \overline{G_{ad}}$ cannot be $\SL_2$-regular.

\subsection{$\SL_2$-regular varieties.}\label{SS:regular}

Recall that an $\SL_2$-regular variety is one which admits an action of $\SL_2$ such that any one-dimensional unipotent subgroup of $\SL_2$ fixes a single point. Aky{\i}ld{\i}z and Carrell developed a remarkable approach for studying the cohomology algebra $H^*(X;\C)$ of such varieties, see~\cite{AC1987}.
Here, we will give a brief summary of how their method works.

Let $\mathbf{B}$ denote the Borel subgroup of upper triangular matrices in $\mathbf{SL}_2$.
Let $\mathbf{U}$ denote its unipotent radical, that is, 
\[
\mathbf{U} := \left\{ \begin{bmatrix} 1 & a \\ 0 & 1 \end{bmatrix} :\ a\in \mathbb{G}_a \right\}.
\] 
Clearly, $\mathbf{U}$ is generated by the regular unipotent element $u:=  \begin{bmatrix} 1 & 1 \\ 0 & 1 \end{bmatrix} $. 
The maximal diagonal torus of $\mathbf{B}$, which is one dimensional, will be denoted by $\mathbf{T}$. 
It is given by the image of the one-parameter subgroup $\lambda : \mathbb{G}_m \to \SL_2$, $t\mapsto \text{diag}(t,t^{-1})$. 
Now let $X$ be a smooth projective $\mathbf{SL}_2$-variety. 
We will denote by $o$ the unique $\mathbf{U}$- (or, equivalently, $u$-fixed) point in $X$. 
Then $\mathbf{T}$ has finitely many fixed points; $o$ is one of these $\mathbf{T}$-fixed points.
Hence, the tangent space at $o$, denoted $T_o X$, has the structure of a $\mathbf{T}$-module. 
All $\mathbf{T}$-weights on $T_o X$ are negative integers.  
Let us label them as follows: $b_1 > b_2 > \cdots > b_k$.  
For $i\in \{1,\dots, k\}$, we set $a_i := - b_i$.
Also, we will denote by $M_{b_i}$ the $\mathbf{T}$-weight subspace corresponding to the weight $b_i$. 
We set $\mu_i := \dim M_{b_i}$, $i\in \{1,\dots, k\}$. 
Aky{\i}ld{\i}z and Carrell show in~\cite{AC1989} that $H^*(X,\C)$ is isomorphic to the finite dimensional graded $\C$-algebra $\C[X_o]/I$,
where $\C[X_o]$ is identified with the symmetric algebra of the cotangent space at $o$, and $I$ is the ideal of the fixed point scheme of $\mathbf{U}$.
It follows that the Poincar\'e polynomial of $X$ is given by 
\begin{align}\label{A:KMformula}
P(X,t) = \prod_{1\leq i \leq k } \left( \frac{1-t^{2(a_i+1)}}{1-t^{2a_i}} \right)^{\mu_i}.
\end{align}
All partial flag varieties are known to be $\SL_2$-regular. 
The Poincar\'e polynomial of a flag variety is given by the formula $P(G/B,t) = \sum_{w\in W} t^{2\ell(w)}$, 
hence, (\ref{A:KMformula}) gives the famous Kostant-Macdonald identity~\cite{AC1989}.

\section{The proof of Theorem~\ref{T:secondthm}}\label{S:second}

Let $G$ be a semi-simple simply connected complex algebraic group.
A unipotent element $u$ is said to be {\em regular} if there is a unique Borel subgroup $B$ of $G$ such that $u\in B$. 
In~\cite[Theorem 1, \S3.1]{Steinberg1974}, Steinberg shows that if $u$ is not a regular unipotent element, then it belongs to infinitely many Borel subgroup. Let us restate this fact in terms of $u$-fixed subvarieties of the flag variety.
\begin{Lemma}\label{L:characterizeregular}
Let $u$ be a unipotent element of $G$. 
Then $(G/B)_u = pt$ if and only if $u$ is regular. 
In this case, the trivial coset $id_G B$ is the Springer fiber $(G/B)_u$. 
\end{Lemma}

\begin{Remark}
It is easy to see that, for every unipotent element $u\in \U(G)$, 
the inclusion of parabolic subgroups $P\subseteq Q$ induces a surjective morphism of closed sets, $(G/P)_u \to (G/Q)_u$.
\end{Remark}

The following result is proved in Borel's book~\cite[Ch IV, Proposition 14.22]{Borel}.
\begin{Lemma}\label{L:P14.22}
Let $G$ be a connected reductive group.
Let $P$ and $Q$ be parabolic subgroups of $G$.
Then the following statements hold:
\begin{enumerate}
\item $P\cap Q$ is connected and $(P\cap Q)\cdot \mathscr{R}_u P$ is a parabolic subgroup.
\item Let $L$ be a Levi subgroup of $P$.
Then $H\mapsto H\cdot \mathscr{R}_u P$ is a bijection of the set of parabolic subgroups of $L$ onto the set of parabolic subgroups of $P$. 
\item If $Q$ is conjugate to $P$ and contains $\mathscr{R}_u P$, then $Q=P$.
\end{enumerate}
\end{Lemma}

\begin{Lemma}\label{L:notin}
Let $u$ be a regular unipotent element contained in a Borel subgroup $B$. 
If $Q$ is a parabolic subgroup of $G$ such that $Q\neq G$, then $u$ is not contained in a Levi factor of $Q$.
\end{Lemma}
\begin{proof}
Towards a contradiction let assume that $u\in L$, where $L$ is a Levi factor of a parabolic subgroup $Q$ such that $Q\neq G$.
Since $u$ is unipotent, it is contained in a maximal unipotent subgroup of $L$.
Hence, $u$ is contained in a Borel subgroup of $L$. 
Let $B_L$ be a Borel subgroup of $L$ that contains $u$. 
It follows from Part 2. of Lemma~\ref{L:P14.22} that $B_L\cdot \mathscr{R}_u Q$ is a Borel subgroup in $Q$, hence in $G$.
In particular, by the uniqueness of the Borel subgroup that contains a regular element, we see that $B= B_L\cdot \mathscr{R}_u Q$.
Since $B_L$ is uniquely determined by $B$, $u$ is a regular unipotent element of $L$. 
Let $T$ be a maximal torus of $B_L$ (hence of $B$).
Let $\Delta$ denote the set of simple roots of $G$ determined by the pair $(B,T)$. 
Let $\Delta_L$ denote the set of simple roots of $L$ determined by the pair $(B_L,T)$.
Then $\Delta_L \subsetneq \Delta$. 
Recall that a character $\alpha$ of $T$ is called a root of $G$ (with respect to $T$) 
if there exists a morphism of algebraic groups $x_\alpha : \mathbb{G}_a \to G$ such that 
the following conditions are satisfied: 
\begin{enumerate}
\item $x_\alpha$ is an isomorphism onto the image $U_\alpha$, which is a unipotent group, normalized by $T$.
\item $t x_\alpha(r) t^{-1} = x_\alpha ( \alpha(t) r)$ for every $t\in T$ and $r\in \C$. 
\end{enumerate}
Now, in this notation, by~\cite[Theorem 1, \S3.7]{Steinberg1974}, $x\in U:=\mathscr{R}_u B$ is a regular unipotent element, 
if $x= \prod_{\alpha > 0} x_\alpha (c_\alpha)$, then $c_\alpha \neq 0$ for every simple root $\alpha$ in $\Delta$. 
But since $U_L:= \mathscr{R}_u B_L$ is a subgroup of $U$, $\Delta_L \subsetneq \Delta$, and the root system of $(L,T)$
is a sub-root system of $(G,T)$, we obtain the equality
\begin{align}\label{A:simpleroots}
 \prod_{\alpha \in \Phi^+} x_\alpha (c_\alpha) =   \prod_{\alpha \in \Phi_L^+} x_\alpha (c_\alpha'),
\end{align}
where $c_\alpha c_\alpha' \neq 0$ for every simple root $\alpha$. 
Since there are more simple roots that appear on the left hand side of (\ref{A:simpleroots}), and since the nonzero values of the morphisms 
$x_\alpha$, $\alpha \in \Phi^+$ are free, the equation (\ref{A:simpleroots}) cannot hold. This contradiction shows that $u$ does not belong to a Levi 
subgroup of $G$ unless $L=G$. 
\end{proof}

\begin{Corollary}\label{C:containedinQ}
Let $u$ be a regular unipotent element of $G$. 
If $Q$ is a parabolic subgroup such that $u \in Q$, then the unique Borel subgroup of $G$ that contains $u$ is contained in $Q$. 
\end{Corollary}

\begin{proof}
Let $B$ be the unique Borel subgroup such that $u\in B$. 
Then $B':=B\cap Q \cdot \mathscr{R}_u Q$ is a Borel subgroup of $Q$ by Lemma~\ref{L:P14.22}.
In particular, $B'$ is a Borel subgroup of $G$. 
Since $u\in B'$, we see that $B' = B$ by the uniqueness of $B$. 
\end{proof}

\begin{Corollary}\label{C:nofixed}
Let $u$ be a regular unipotent element in $G$.
Let $Q$ be a parabolic subgroup such that $Q\neq G$. 
If $L$ is a Levi factor of $Q$, then $(G/L)_u = \emptyset$. 
\end{Corollary}

\begin{proof}
Let $g$ be any element from $G$. 
If $B$ denotes the unique Borel subgroup of $G$ such that $u\in B$, then $gBg^{-1}$ is the unique Borel subgroup of $G$ that contains $gug^{-1}$. 
In particular, $gug^{-1}$ is a regular unipotent element of $G$. 
We now assume that there is an element $gL$ in $(G/L)_u$.
Then $g^{-1}ug$ is a regular unipotent element of $G$, and it is contained in $L$. 
But this contradicts with Lemma~\ref{L:notin}, hence, the proof of our assertion is complete. 
\end{proof}

Two parabolic subgroups of $G$ are said to be {\em opposite} if their intersection is a common Levi subgroup.
\begin{Lemma}\label{L:P14.21}
Let $P$ be a parabolic subgroup of $G$.
If $L$ is a Levi subgroup of $P$, then there exists one and only one parabolic subgroup $P'$ in $G$ such that $P'\cap P= L$. 
Any two parabolic subgroups opposite to $P$ are conjugate by a unique element of $\mathscr{R}_u P$. 
\end{Lemma}
\begin{proof}
This is proved in~\cite[Ch IV, Proposition 14.21]{Borel}.
\end{proof}

\begin{Theorem}\label{T:finitenormalizer}
Let $G$ be reductive group, not necessarily connected.
Let $P$ be a parabolic subgroup of $G$, and let $L$ be a Levi factor of $P$. 
Then the quotient group $N_G(L) / L$ is finite. 
\end{Theorem}

\begin{proof}
Let $L$ be a Levi subgroup in any reductive group $G$. 
Let $S$ be the central torus of $L$. 
Then the centralizer $Z_G(S)$ of $S$ in $G$ equals $L$.  
Also, any element $g$ of $N_G(L)$ normalizes the center of $L$ and hence it normalizes the central torus $S$. 
Conversely, any element $g$ that normalizes $S$ will normalize the centralizer $L$ of $S$. 
Thus, we conclude that $N_G(S) = N_G(L)$. 
Now let $T$ be a maximal torus of $G$ containing $S$. 
As $T$ commutes with $S$, it is contained in $L$. 
Let $g$ be an element of $N_G(S) = N_G(L)$. 
Then the torus $gTg^{-1}$ contains $S$ and hence it commutes with $S$. 
Therefore, $gTg^{-1}$ is contained in $Z_G(S) = L$. 
Thus, both $T$ and $gTg^{-1}$ are maximal tori in $L$. 
By conjugacy of maximal tori in $L$, we see that there is an element $x$ in $L$ such that  $T = xgTg^{-1}x^{-1}$.
In other words, the element $xg$ is contained in $N_G(T)$. 
Thus, we see that $g$ belongs to $x^{-1}N_G(T) \subset LN_G(T)$. 
Therefore, we have $N_G(L)\subseteq LN_G(T)$. 
Since the Weyl group $W = N_G(T)/T$ is finite and $T$ is contained in $L$, we infer now that $N_G(L)/L$ is finite.
\end{proof}

A spherical subgroup $H\subseteq G$ is said to be {\em sober} if $N_G(H)/H$ is finite (or, equivalently, the ``valuation cone'' $\mathcal{V}(G/H)$ is strictly convex).  
\begin{Corollary}\label{C:finitenormalizer}
Let $G$ be a reductive group.
Let $P$ be a parabolic subgroup of $G$, and let $L$ be a Levi factor of $P$. 
If $L$ is a spherical subgroup of $G$, then $L$ is a sober subgroup. 
\end{Corollary}
\begin{proof}
This follows from Theorem~\ref{T:finitenormalizer}.
\end{proof}

\begin{Remark}
The concept of sobriety is crucial for our purposes since it is a necessary condition for the existence of the wonderful completions. 
If $H$ is a sober subgroup, then the {\em standard embedding of $G/H$}, which we denote by $O$, is defined as the unique toroidal simple complete $G$-equivariant embedding $G/H \hookrightarrow O$ defined by the ``colored cone'' $(\mathcal{V}(G/H),\emptyset)$. 
The {\em wonderful completions} are precisely the smooth standard embeddings. 
\end{Remark}

We are now ready to prove our second main result. Let us recall its statement.

Let $G$ be a connected reductive group, and let $L$ be a Levi subgroup of a parabolic subgroup $Q\subsetneq G$. 
If $u$ is a regular unipotent element of $G$, then $u$ acts freely on $G/J$, where $J$ is any closed subgroup of the normalizer of $L$ in $G$.

\begin{proof}[Proof of Theorem~\ref{T:secondthm}]
First, we will prove that $(G/N_G(L))_u=\emptyset$. 
We assume towards a contradiction that $(G/N_G(L))_u$ is nonempty. 
Since $N_G(L)/L$ is finite (Theorem~\ref{T:finitenormalizer}), the canonical projection $\pi: G/L \to G/N_G(L)$ is a principal $N_G(L)/L$-bundle with finite fibers. In particular, it is a finite morphism. 
Therefore, it maps closed sets to closed sets. 
It follows that $\pi ( (G/N_G(L))_u )$ is a nonempty closed subset of $G/L$.
But $\pi$ is $G$-equivariant, therefore, any element of $\pi ( (G/N_G(L))_u )$ is fixed by $u$. 
This contradicts with the conclusion of Corollary~\ref{C:nofixed}.

Now, to finish the proof, let $J$ be a closed subset of $N_G(L)$. 
Then we have an equivariant surjective morphism $G/J\to G/N_G(L)$. 
If $u$ fixes a point in $G/J$, then it does fix a point in $G/N_G(L)$.
Since  $(G/N_G(L))_u = \emptyset$, the proof of our assertion follows. 
\end{proof}

\section{The proof of Theorem~\ref{T:firstthm}}\label{S:first}

First, let us recall our Theorem~\ref{T:firstthm}: If $L$ is a Levi subgroup of the form $L=G^\sigma$ for some involution $\sigma$ of $G$, then the wonderful completion of $G/N_G(L)$ is an $\SL_2$-regular variety. 


\begin{proof}[Proof of Theorem~\ref{T:firstthm}]
Let $X$ denote the wonderful completion of $G/N_G(L)$. 
Let $T$ denote the maximal $\sigma$-stable torus in $G$ such that $T$ contains a maximal $\sigma$-anisotropic torus.
Every $G$-orbit in $X$ is fibration over an appropriate partial flag variety $G/Q$. 
Let $L_Q$ denote the Levi factor of $Q$ such that $T\subseteq L_Q$. 
If $M$ denotes the derived subgroup of the Levi factor $L_Q$ of $Q$, then the $G$-orbit $O$ that corresponds to $Q$ is given by 
$O \cong G\times_Q (M/M^\sigma)$, where $\sigma$ is the induced automorphism on $M$. 
Then $\sigma |_M$ is either a nontrivial involution, or it is the identity automorphism.
In latter case, $M/M^\sigma$ is just a point, hence, the $G$-orbit $O$ is the closed orbit $G/Q$. 
If $\sigma|_M$ is a nontrivial involution, then since $L$ is a Levi subgroup of $G$, the subgroup $M^\sigma$ is a Levi subgroup of $M$ as well. 
Indeed, $M^\sigma$ is the intersection of $L$ with a $\sigma$-stable Levi subgroup of $Q$, see~\cite[Section 5.2]{DP}.
In particular, by Corollary~\ref{C:nofixed}, the induced unipotent transformation on $M/M^\sigma$, has no fixed points. 
It follows that $u$ has no fixed points in $O$. 
This means that the only fixed points of $u$ are on the closed orbit.
Therefore, the set of $u$-fixed points of $X$ is equal to the set of $u$-fixed points in a partial flag variety.
As we have mentioned before, this fixed subset is a singleton. 
Finally, since the unipotent subgroup generated by $u$, denoted $\langle u \rangle$, is normalized by $T$, the Jacobson-Morozov Theorem 
guarantees the existence of a copy of $\SL_2$ in $G$ containing $\langle u \rangle$ as a maximal unipotent subgroup. 
Hence, we have a regular $\SL_2$-action on $X$. 
This finishes the proof of our theorem. 
\end{proof}

Let $X$ be a wonderful completion as in Theorem~\ref{T:firstthm}. 
Let $Q$ denote the stabilizer of a point in the unique closed $G$-orbit. 
Then $G/Q$ is a flag variety.  
Following our notation from Subsection~\ref{SS:regular}, let us denote by $o$ the unique $\mathbf{U}$-fixed point. 
This point is the unique $u$-fixed point $id_GQ$ in $G/Q$. 
As in the proof of Theorem~\ref{T:firstthm}, we view $\mathbf{U}$ as a maximal unipotent subgroup of a copy of $\SL_2$ in $G$;
the maximal torus $T$ of $L$ normalizes $\mathbf{U}$, hence, the image of the one-parameter group $\lambda: \mathbb{G}_a \to \mathbf{T}$ 
is contained in $T$.

We are now ready to prove Corollary~\ref{C:firstcorollary}, which states that the Poincar\'e polynomial of $X$ is divisible by the Poincar\'e polynomial of the closed orbit. 

\begin{proof}[Proof of Corollary~\ref{C:firstcorollary}]

To compute the Poincar\'e polynomial, first, we will determine the $\mathbf{T}$-weights at the tangent space $T_o X$. 
The general procedure for this calculation is presented in the articles~\cite{DP,DS}.  
Since we are working with the torus fixed point in the closed orbit, in our situation the calculation is simpler. 
In the notation of~\cite{DS}, the tangent space of $X$ at $o$ has a $T$-equivariant (hence $\mathbf{T}$-equivariant) decomposition 
into three subspaces, $T_o X= T_1 \oplus T_2 \oplus T_3$.
The first two summands of this decomposition come from the fibration of a $G$-orbit over a partial flag variety. 
In our case, the fibration map is the identity map, therefore, $T_1$ is the tangent space at $o$, that is $T_o (G/Q)$.
The second summand is the tangent space in the ``vertical direction'', hence, it is trivial, $T_2 = \{0\}$.
Finally, $T_3$ is isomorphic to the normal space $T_o X / T_o (G/Q)$. 
Now from Subsection~\ref{SS:regular}, we know that the Poincar\'e polynomial $P(X,t)$ is given by a product over all $\mathbf{T}$-weights on the tangent space $T_o X$. 
As a result of the equivariant decomposition, we see that $P(X,t)$ has a factor that comes from the $\mathbf{T}$-module structure of the tangent space $T_o (G/Q)$.
But according to~\cite[Corollary 1 and remarks following it]{AC2012}, this is precisely the Poincar\'e polynomial of the closed orbit $G/Q$. 
At the same time, we know from Subsection~\ref{SS:wonderful embeddings} that the spherical weights $\varSigma(X)$ are precisely the $T$-weights 
on the normal space $T_o X / T_o (G/Q)$. 
The Local Structure Theorem~\cite{BLV1987} adapted to wonderful completions imply that the dimension of $X$ is equal to $\dim G/Q + \dim (T_o X / T_o (G/Q))$. 
Since every $T$-weight is a $\mathbf{T}$-weight, and since all of the spherical weights are distinct, we see that the $\mathbf{T}$-weights on $T_o X / T_o (G/Q)$ are all distinct. Hence, it is the (unique) irreducible $\SL_2$-module of dimension $r:=|\varSigma(X)|$. It follows that $a_1=1,\dots , a_r=r$, and $\mu_i =1$ for all $i\in \{1,\dots,r\}$. 
Then, by using the Aky{\i}ld{\i}z and Carrell formula (\ref{A:KMformula}), we see that the contribution of the normal space at $o$ to the Poincar\'e polynomial is $\prod_{ i=1}^r \frac{1-t^{2(i+1)}}{1-t^{2i}}$.
This finishes the proof of our corollary.
\end{proof}

\bibliographystyle{plain}
\bibliography{references}

\end{document}